\documentclass[12pt]{article}
\usepackage{amsxtra,amssymb,amsthm,amsmath,latexsym}

\theoremstyle{plain}

\def\R{{\mathbb R}}

\def\oH{{\overset{\circ}{H}}}
\def\oH1{{\overset{\circ}{H}\kern-.02in{}^1}}

\def\bee{\begin{equation*}}
\def\eee{\end{equation*}}
\def\be{\begin{equation}}
\def\ee{\end{equation}}

\begin{document}

\title{Symmetry problems in harmonic analysis 
}

\author{Alexander G. Ramm\\
 Mathematics Department, Kansas State University, \\
 Manhattan, KS 66506-2602, USA\\
ramm@ksu.edu\\
http://www.math.ksu.edu/\,$\sim$\,ramm}

\date{}
\maketitle\thispagestyle{empty}

\begin{abstract}

Symmetry problems in harmonic analysis are formulated and solved. One of these problems is equivalent to the refined Schiffer's conjecture which was recently proved by the author.

Let $k=const>0$ be fixed, $S^2$ be the unit sphere in $\R^3$, $D$ be a connected bounded domain with $C^2-$smooth boundary $S$, $j_0(r)$ be the spherical Bessel function.

The harmonic analysis symmetry problems are stated in the following theorems. 

{\bf Theorem A.} {\em Assume that $\int_S e^{ik\beta \cdot s}ds=0$ for all $\beta\in S^2$. Then $S$ is a sphere of radius $a$, where $j_0(ka)=0$. } 

{\bf Theorem B.} {\em Assume that $\int_D e^{ik\beta \cdot x}dx=0$ for all $\beta\in S^2$.
 Then $D$ is a ball.}

\end{abstract}

\section{Introduction}\label{S:1}

Symmetry problems for PDE were studied in many
publications by many authors, see,  \cite{R670}, \cite{R691}
and references therein.

Throughout we assume that $D$ is a bounded connected 
$C^2-$smooth domain in $\R^3$, $S$ is the boundary of $D$, 
$N$ is the unit normal to $S$, pointing out of $D$, $u_N$ is the normal derivative of $u$ on $S$, $D'=\R^3\setminus D$, $S^2$ is the unit sphere 
in $\R^3$,  $k>0$ and $c$ are fixed constants, $j_{0}(r)$ is the spherical Bessel function 
and  $g(x,y,k):=\frac{e^{ik|x-y|}}{4\pi |x-y|}$. By 
$(u_N)_{-}$ we denote the limiting value on $S$ of 
$u_N$ from $D'$ and by $(u_N)_{+}$ we denote the limiting value on $S$ of 
$u_N$ from $D$.    

In \cite{R691}, \cite{R694} the refined Schiffer's conjecture (SC) is proved. Let us formulate this result.

{\bf Theorem 1.} {\em Assume that
\be\label{e1}
\Delta u+k^2u=0 \quad in \quad D, \quad u|_S=0,\quad u_N=c.
\ee
Then  $S$ is a sphere of radius $a$ such that 
 $j_0(ka)=0$.}
 
Let us formulate our new result: formulation and solution of a symmetry problem  in harmonic analysis  (Problem HA):
 
 {\bf Theorem A.}{\em Assume that
\be\label{e2}
\int_Se^{ik\beta \cdot y}dy=0\quad \forall \beta \in  S^2.
\ee
Then  $S$ is a sphere of radius $a$, where $j_0(ka)=0$. } 
 
{\bf Theorem B.} {\em Assume that 
\be\label{e3}
\int_D e^{ik\beta \cdot x}dx=0\quad \forall \beta\in S^2,
\ee 
where $D$ is a bounded connected domain in $\R^3$ and  $S^2$ is the unit sphere in $\R^3$. Then $D$ is a ball.}

 We prove that the harmonic analysis symmetry problem (HA),
 Theorem A, is equivalent to the refined Schiffer's conjecture (SC), Theorem 1: if Theorem A holds then Theorem 1 holds and vice versa.

The author does not know any symmetry results in harmonic analysis of the type presented in Theorems A and B.

Theorem A  says  that if the Fourier transform of a distribution supported  on a smooth closed surface $S$
with a constant density has a spherical surface of zeros, then $S$ is a sphere.

Theorem B says that if the Fourier transform of a
characteristic function of a connected bounded domain $D$ has a spherical surface of zeros, then $D$ is a ball. 
 
In Section 2 proofs are given.

\section{Proofs}\label{S:2}

If problem \eqref{e1} has a solution then this solution is unique by the uniqueness of the solution to the Cauchy problem for the Helmholtz elliptic equation \eqref{e1}.  

The solution to equation \eqref{e1} by the Green's formula is:
\be\label{e4}
u(x)=c\int_{S}g(x,t)dt, \quad x\in D; \quad  u(x):=c\int_{S}g(x,t)dt=0, \quad x\in D'.
\ee
Formulas \eqref{e4} are obtained by the 
standard application of the Green's formula.

Namely, one starts with the equations
\be\label{e41}
(\nabla^2_y +k^2)u=0,
\ee
\be\label{e42}
(\nabla^2_y +k^2)g(x,y)=-\delta(x-y), \quad y\in D.
\ee
Multiply \eqref{e41} by $g=g(x,y)$, equation \eqref{e42}
by $u(y)$, subtract the second equation from the first, integrate over $D$ and use the definition of the
delta-function $\delta(x-y)$ and the boundary conditions 
in \eqref{e1} to get \eqref{e4}.

The function $u$, defined by the first formula \eqref{e4}
in $\R^3$ satisfies the radiation condition
\be\label{e43}
 u_r-iku=O(|x|^{-2}), \quad r:=|x|\to \infty
 \ee
 uniformly with respect to directions of $x$.
 
Let $B_R=\{x: |x|\le R\}$, $D\subset B_R$.
If $D$ is a ball $B_a$ of radius $a$, and $u$ solves 
\eqref{e1} then $a$ solves the equation $j_0(ka)=0$,
and the solution $u$ has the form:
\be\label{e5}
u=c\frac{j_0(kr)}{kj_0'(ka)}, \quad r=|x|,
\ee
where $j'_0(r):=\frac {dj_0(r)}{dr}$.

{\bf Proof of Theorem A.} 

Assume that \eqref{e2} holds. Let $u$ be defined by
the first formula \eqref{e4} in $\R^3$. Then, due to \eqref{e2}, one has: 
\be\label{e6}
u=O(|x|^{-2}), \quad |x|\to \infty,
\ee
and \eqref{e43} holds.
Moreover, $u$, defined in \eqref{e4}, solves the equation
\be\label{e7}
(\nabla^2 +k^2)u=0 \quad in \quad D'.
\ee
By the known lemma, see, for example, \cite{R670}, p.30, Lemma 1.2.1, it follows from \eqref{e6}, \eqref{e43} and \eqref{e7} that $u=0$ in $D'$.

For convenience of the reader let us formulate the lemma we
have used.

{\bf Lemma 1.} {\em If \eqref{e6} and \eqref{e7}
hold, then $u=0$ in $D'$.}

Since $u=0$ in $D'$, $u$ is a single layer potential and $S$ is $C^2-$smooth, one concludes that $u$ is continuous
up to $S$ together with its first derivatives, so
\be\label{e8}
u=0, \quad (u_N)_{-}=0 \quad on \quad S.
\ee
By the jump formula for the normal derivative of $u$
(see, for example, \cite{R670}, p. 18), one gets:
\be\label{e9}
(u_N)_{+}-(u_N)_{-}=(u_N)_{+}=1,
\ee
since the density of the single layer potential $u$ is equal to $1$ and $(u_N)_{-}=0$.

Therefore, if \eqref{e2} holds, then $u$ solves problem \eqref{e1} with $c=1$. Consequently, by Theorem 1, $S$ is a sphere
of radius $a$, where $j_0(ka)=0$.

Theorem A is proved.   \hfill$\Box$

{\bf Lemma 2.} {\em Theorem 1 and Theorem A are equivalent.} 

In the proof of Lemma 2 we use the following formula:
\be\label{e10}
g(x,y,k)=\frac {e^{ik|x|}}{4\pi |x|}e^{ik\beta \cdot y}
 +O\big(\frac 1 {|x|^2}\big), \quad |x|\to \infty, \quad \beta=-x/|x|,
\ee
where $|y|\le R$. 


{\em Proof of Lemma 2.} Assume that Theorem 1 holds. 
  Define $u$ by formula \eqref{e4}. As $|x|\to \infty$, $x/|x|=-\beta$, this yields \eqref{e2}. So, if
Theorem 1 holds, then Theorem A holds.

Conversely, Suppose that Theorem A holds. From \eqref{e2}
one derives the relation:
\be\label{e11}
u(x):=\int_S g(x,t)dt=0\quad in \quad D'.
\ee
Indeed, the integral $u(x)$ in \eqref{e11}
satisfies differential equation \eqref{e7}
in $D'$ and
 $u=O(|x|^{-2})$ as $|x|\to \infty$.
So $u=0$ in $D'$ by Lemma 1.
Equation \eqref{e7} for $u$ holds in $D$, $u=0$ on $S$
by continuity, and $(u_N)_{+}=1$ on $S$ by the jump formula for the normal derivatives of the single layer potential $u$.
Thus, $u$ solves problem \eqref{e1}. So, Theorem A yields
the conclusion of Theorem 1.

 Lemma 1 is proved. \hfill$\Box$

{\bf Proof of Theorem B.} Assume that \eqref{e3} holds.
Define
\be\label{e12}
w(x):=\int_D g(x,t)dt, \quad x\in \R^3.
\ee
Then
\be\label{e13}
(\nabla^2+k^2)w=0 \quad in \quad D',
\ee
and 
\be\label{e14}
w=O(|x|^{-2}) \quad as\quad |x|\to \infty.
\ee
Therefore, by Lemma 1, one concludes that 
\be\label{e15}
w=0 \quad in \quad D'.
\ee
Since $w$ is a volume potential which is continuous together with its first derivatives in $\R^3$,
one gets from \eqref{e15} and \eqref{e12} that
\be\label{e16}
w=0, \quad w_N=0\quad on \quad S,
\ee  
and 
\be\label{e17}
(\nabla^2+k^2)w=-1 \quad in \quad D.
\ee
We now use Theorem 3.1 from \cite{R691}, p.15, and conclude that $D$ is a ball.

Theorem B is proved. \hfill$\Box$

For convenience of the reader let us formulate Theorem 3.1 from \cite{R691}. The assumptions about $D$ are the same
as in this paper. Below $c_j$, $j=0,1,2,$ are some constants.

{\bf Theorem 3.1.} {\em Assume that the problem
\be\label{e18}
(\nabla^2+k^2)w=c_0 \quad in\quad D, \quad u|_S=c_1,\quad u_N|_S=c_2,
\ee
is solvable. If
\be\label{e19}
|c_1-c_0k^{-2}|+|c_2|>0,
\ee
then $D$ is a ball.}

In our case $c_1=c_2=0$ and $c_0=-1$, so condition \eqref{e19} is satisfied.

\newpage

\end{document}